\documentclass[3p,onecolumn]{article}
\usepackage[utf8]{inputenc}
\usepackage[T1]{fontenc}
\usepackage{lmodern}
\usepackage{amsmath}
\usepackage{amsfonts}
\usepackage{amssymb}
\usepackage{amsthm}
\usepackage[english]{babel}
\newcommand{\bm}[1]{\boldsymbol{#1}}

\newcommand{\R}{\mathbb{R}}
\newcommand{\A}{\mathcal{A}}
\newcommand{\F}{\mathcal{F}}
\newcommand{\Pnm}{\mathcal{P}(n,m)}
\newcommand{\ulf}{\underline{f}}
\newcommand{\olf}{\overline{f}}
\newtheorem{theorem}{Theorem}
\newtheorem{lemma}{Lemma}
\newtheorem{corollary}{Corollary}

\begin{document}

\title{Boundary optimization for rough sets}

\author{Konrad Engel \thanks{Universit\"at Rostock,  Institut f\"ur Mathematik, 18051 Rostock, Germany. E-mail: konrad.engel@uni-rostock.de} 
\and
Tran Dan Thu\thanks{School of Information Technology, University of Science Ho Chi Minh City, 227 Nguyen Van Cu, Dist. 5, Ho Chi Minh City, Vietnam. 
E-mail: tdt@hcmus.edu.vn}
}

\maketitle

\begin{abstract}
Let $n > m\ge 2$ be integers and let $\A=\{A_1,\dots,A_m\}$ be a partition of $[n]=\{1,\dots,n\}$.
For $X \subseteq [n]$, its $\A$-boundary region $\A(X)$ is defined to be the union of those blocks $A_i$ of $\A$ for which $A_i\cap X\neq \emptyset$ and 
$A_i\cap ([n] \setminus X)\neq \emptyset$. For three different probability distributions on the power set of $[n]$, partitions $\A$ of $[n]$ are determined such that the expected cardinality of the $\A$-boundary region of a randomly chosen subset of $[n]$ is minimal and maximal, respectively. The problem can be reduced to an optimization problem for integer partitions of $n$. In the most difficult case, the concave-convex shape of the corresponding weight function as well as several other inequalities are proved using an integral representation of the weight function. In one case, there is an interesting analogon to the AZ-identity.
The study is motivated by the rough set theory.
\end{abstract}


\section{Introduction}
Let $n > m\ge 2$ be integers and let $\A=\{A_1,\dots,A_m\}$ be a partition of $[n]=\{1,\dots,n\}$.
For $X \subseteq [n]$, the \emph{$\A$-lower approximation of $X$} is defined by
\[
\A^-(X)=\bigcup_{i: A_i \subseteq X} A_i
\]
and the \emph{$\A$-upper approximation of $X$} is defined by
\[
\A^+(X)=\bigcup_{i: A_i \cap X\neq \emptyset} A_i\,.
\]
Obviously,  $\A^-(X) \subseteq X \subseteq \A^+(X)$ for all $X$.
The set
\[
\A(X)=\A^+(X) \setminus \A^-(X)
\]
is called the \emph{$\A$-boundary region of $X$}.

These notions stem from \emph{rough set theory} that was initiated by Pawlak in \cite{Paw82}, see also \cite{Paw91}.
Here $[n]$ is as a set of objects (the universe) and, in addition, the values of a given set $P$ of attributes are known for each object. Then the partition $\A$ of $[n]$ is the set of classes of that equivalence relation, where two objects are related if both objects have the same attribute-values.
Now assume that there is given a set $X$ of objects, but we have only access to the attribute-values of the objects in $X$. Then $\A^-(X)$ consists of those objects that can be positively identified as members of $X$ and $\A^+(X)$ consists of those objects that might be members of $X$.
Thus $|\A(X)|$ can be interpreted as a measure of uncertainty of $X$ generated by the attribute-values. If $\A(X) = \emptyset$ then $X$ is called a \emph{crisp set}, otherwise it is called a \emph{rough set}, but we consider crisp sets also as degenerated rough sets.

In this paper, we provide sharp bounds on the average size of the $\A$-boundary region with respect to three reasonable distributions.
\begin{enumerate}
\item Uniform distribution, i.e., each set $X$ has probability $P_1(X) = \frac{1}{N_1}$, where $N_1 = \sum_{k=0}^n \binom{n}{k} = 2^n$. Let
\[
\mu_1(\A)=\frac{1}{N_1} \sum_{X \subseteq [n]} |\A(X)|\,.
\]
This number can be interpreted as an average \emph{absolute uncertainty}.
\item ``Relative distribution'', i.e., each set $X\neq \emptyset$ has probability $P_2(X) = \frac{1}{N_2}\frac{1}{|X|}$, where $N_2 = \sum_{k=1}^n \binom{n}{k} \frac{1}{k}$. (For $X= \emptyset$, we set $P_2(X) =0$). Let
\[
\mu_2(\A)=\frac{1}{N_2} \sum_{\emptyset \neq X \subseteq [n]} \frac{1}{|X|}|\A(X)|\,.
\]
This number can be interpreted as an average \emph{relative uncertainty}.
\item ``Size dependent relative distribution'', i.e., each set $X\neq \emptyset$ has probability $P_3(X) = \frac{1}{N_3}\frac{1}{\binom{n}{|X|}}\frac{1}{|X|}$, where $N_3 = \sum_{k=1}^n \frac{1}{k}$. (For $X= \emptyset$, we set $P_3(X) =0$). Let
\[
\mu_3(\A)=\frac{1}{N_3} \sum_{\emptyset \neq X \subseteq [n]} \frac{1}{\binom{n}{|X|}|X|}|\A(X)|\,.
\]
This number can be interpreted as the average of the \emph{average relative uncertainty for fixed cardinality of $X$}.
\end{enumerate}
Now the question is, how small and how large the numbers $\mu_j(\A)$ can be if all partitions of $[n]$ into $m$ parts are allowed
(throughout let $j \in \{1,2,3\}$).

\section{Reduction to an optimization problem for integer partitions}

Since the factors $N_j$ do not have influence, we define
\begin{align*}
f_1(\A)&= \sum_{X \subseteq [n]} |\A(X)|\,,\\
f_2(\A)&=\sum_{\emptyset \neq X \subseteq [n]} \frac{1}{|X|}|\A(X)|\,,\\
f_3(\A)&=\sum_{\emptyset \neq X \subseteq [n]} \frac{1}{\binom{n}{|X|}|X|}|\A(X)|\,
\end{align*}
and further
\begin{align*}
\ulf_j&=\min\{f_j(\A): \A \text{ is a partition of } [n] \text{ into $m$ parts}\}\,,\\
\olf_j&=\max\{f_j(\A): \A \text{ is a partition of } [n] \text{ into $m$ parts}\}\,.
\end{align*}
Partitions $\A$ of $[n]$ with $f_j(\A)= \ulf_j$ or $f_j(\A)= \olf_j$ are called \emph{$j$-minimal} and \emph{$j$-maximal}, respectively.

With a partition $\A=\{A_1,\dots,A_m\}$ we associate its block sizes $a_i=|A_i|, i=1,\dots,m$.
It is not difficult to see that the values $f_j(\A)$ only depend on the block sizes of $\A$ and not on the concrete partition. Indeed, let $w:\{0,\dots,n\} \rightarrow \R$ be any (weight) function, let
\[
f_w(\A)=\sum_{X \subseteq [n]} w(|X|) |\A(X)|
\]
and let $F_w: \{0,\dots,n\} \rightarrow \R$ be defined (with $\binom{n}{k}=0$ if $k < 0$) by
\[
F_w(a)=\sum_{k=0}^n \left( \binom{n-a}{k}+\binom{n-a}{k-a}\right) w(k).
\]
Now the following lemma converts the boundary optimization problem for rough sets into an optimization problem for integer partitions.
\begin{lemma}
\label{lemma1}
We have
\[
f_w(\A)=n \sum_{k=0}^n \binom{n}{k} w(k) - \sum_{i=1}^m a_iF_w(a_i).
\]
\end{lemma}
\proof
We  have
\begin{align*}
f_w(\A)&=\sum_{X \subseteq [n]} w(|X|) \left(|\A^+(X)|-|\A^-(X)|\right)\\
			 &=\sum_{X \subseteq [n]} w(|X|) \left(\sum_{i: A_i \cap X \neq \emptyset} |A_i| - \sum_{i: A_i \subseteq X} |A_i| \right)\\
			 &=\sum_{X \subseteq [n]} w(|X|) \left(n-\sum_{i: A_i \cap X = \emptyset} |A_i| - \sum_{i: A_i \subseteq X} |A_i| \right)\\
			 &=n \sum_{k=0}^n \binom{n}{k} w(k) -\sum_{i=1}^m a_i \left(\sum_{X: X \subseteq [n]\setminus A_i} w(|X|) + \sum_{X: A_i \subseteq X} w(|X|) \right)\\
			 &=n \sum_{k=0}^n \binom{n}{k} w(k) -\sum_{i=1}^m a_i \left(\sum_{k=0}^n \left(\binom{n-a_i}{k} + \binom{n-a_i}{k-a_i}\right) w(k)\right)\\
			 &=n \sum_{k=0}^n \binom{n}{k} w(k) - \sum_{i=1}^m a_iF_w(a_i)\,.
\end{align*}
\qed

 Let $\Pnm$ be the set  of all integer partitions $\bm{a}=(a_1,a_2,\dots,a_m)$ of the integer $n$ into $m$ parts, i.e., 
$a_1,a_2,\dots,a_m$ are  positive integers with 
$a_1+\dots+a_m=n$. With the concrete weight functions $w_j: \{0,...,n\} \rightarrow \R$, where
\begin{align*}
w_1(k)&=1\qquad\text{for all } k\,,\\
w_2(k)&=
\begin{cases}
0, &\text{ if } k=0\,,\\
\frac{1}{k},& \text{ otherwise}\,,
\end{cases}\\
w_3(k)&=
\begin{cases}
0, &\text{ if } k=0\,,\\
\frac{1}{\binom{n}{k}k},& \text{ otherwise}\,,
\end{cases}
\end{align*}
we obtain from Lemma \ref{lemma1}
\begin{align*}
\ulf_j&=n \sum_{k=0}^n \binom{n}{k} w_j(k) - \max\left\{ \sum_{i=1}^m a_i F_{w_j}(a_i): \bm{a} \in \Pnm\right\}\,,\\
\olf_j&=n \sum_{k=0}^n \binom{n}{k} w_j(k) - \min\left\{ \sum_{i=1}^m a_i F_{w_j}(a_i): \bm{a} \in \Pnm\right\}\,.
\end{align*}
We call an integer partition $\bm{a} \in \Pnm$ \emph{$j$-minimal} and \emph{$j$-maximal} if $ \sum_{i=1}^m a_i F_{w_j}(a_i)$ attains the minimum and maximum, respectively.

\begin{corollary}
The set partition $\A$ is $j$-minimal (resp. $j$-maximal) iff the block sizes of $\A$ form a $j$-maximal (resp. $j$-minimal) integer partition $\bm{a} \in \Pnm$.
\end{corollary}
Thus we restrict ourselves to the determination of $j$-optimal integer partitions.

We mention that, for arbitrary weight functions $F$ (instead of the special weight functions $F_{w_j}$) and arbitrary number of items (instead of fixed number of items), a dynamic programming algorithm for the determination of optimal integer partitions is given in \cite{ERS14}. Moreover, \cite{ERS14} contains bounds for the number of pairwise different items in an optimal integer partition. But for our special cases of functions $F_{w_j}$ we may provide explicit solutions.
For more information on integer partitions see \cite{And98}.

In the following we present integer partitions from $\Pnm$ also in the form
$1^{\lambda_1} 2^{\lambda_2} \dots n^{\lambda_n}$ which means that there are exactly $\lambda_i$  items $i$, $i=1,\dots,n$, and, equivalently,
\[
\sum_{i=1}^n i \lambda_i=n \text{ and } \sum_{i=1}^n \lambda_i=m \text{ with non-negative integers $\lambda_1,\dots,\lambda_n$}.
\]
Moreover, items of the form $i^0$ will be mostly omitted.

\section{Main results}

The main results are given by the following three theorems.
\begin{theorem}
\label{theorem1}~
\renewcommand{\labelenumi}{\alph{enumi})}
\begin{enumerate}
\item An integer partition of $\Pnm$ is 1-minimal if it is of the form
\[
\begin{cases}
1^{m-1} (n-m+1)^1&\text{ if } n-m \le 5\,,\\
q^{m-r} (q+1)^r&\text{ if } n-m > 5, n >5m, n=qm+r, 0 \le r <m\,,\\
1^{m-\frac{n-m}{4}} 5^{\frac{n-m}{4}}&\text{ if }n-m > 5, n \le 5m \text{ and }m-n \equiv 0(\hspace{-2ex}\mod 4)\,,\\
1^{m-\frac{n-m+1}{4}} 4^1 5^{\frac{n-m-3}{4}}&\text{ if }n-m > 5, n \le 5m \text{ and }m-n \equiv 1(\hspace{-2ex}\mod 4)\,,\\
1^{m-\frac{n-m+2}{4}} 4^2 5^{\frac{n-m-6}{4}}&\text{ if }n-m > 5, n \le 5m \text{ and }m-n \equiv 2(\hspace{-2ex}\mod 4)\,,\\
1^{m-\frac{n-m-1}{4}} 5^{\frac{n-m-5}{4}} 6^1&\text{ if }n-m > 5, n \le 5m \text{ and }m-n \equiv 3(\hspace{-2ex}\mod 4)\,.\\
\end{cases}
\]
\item An integer partition of $\Pnm$ is 1-maximal if it is of the form
\[
\begin{cases}
2^{m-1} (n-2m+2)^1&\text{ if } n \ge 2m,\\
1^{2m-n} 2^{n-m}&\text{ if } n< 2m\,.
\end{cases}
\]
\end{enumerate}
\end{theorem}

\begin{theorem}
\label{theorem2}
\renewcommand{\labelenumi}{\alph{enumi})}~
Let $n \ge 13$.
\begin{enumerate}
\item An integer partition of $\Pnm$ is 2-minimal if it is of the form
given in part a) of Theorem \ref{theorem1}.

\item An integer partition of $\Pnm$ is 2-maximal if it is of the form
\[
\begin{cases}
2^{m-1}(n-2m+2)^1 &\text{ if } n > 2m+2\,,\\
2^{m-2}3^2&\text{ if } n=2m+2\,,\\
2^{m-1}3^1&\text{ if } n=2m+1\,,\\
1^{2m-n} 2^{n-m}&\text{ if } n \le 2m\,.
\end{cases}
\]
\end{enumerate}
\end{theorem}

\begin{theorem}
\label{theorem3}
\renewcommand{\labelenumi}{\alph{enumi})}~
\begin{enumerate}
\item An integer partition of $\Pnm$ is 3-minimal if it is of the form
\[
1^{m-1} (n-m+1)^1\,.
\]
\item An integer partition of $\Pnm$ is 3-maximal if it is of the form
\[
q^{m-r} (q+1)^r, \text{ where }n=qm+r, 0 \le r <m.
\]
\end{enumerate}
\end{theorem}
\section{Optimization for integer partitions}
It will turn out that our functions $F_{w_j}$ have a concave-convex shape. Therefore we present an elementary optimality criterion that is at least for the special case $c=1$ folklore. Let $g: [n] \rightarrow \R$ be any function. The integer partition $\bm{a} \in \Pnm$ is called \emph{$g$-minimal} (resp. \emph{$g$-maximal}) if the objective function $\sum_{i=1}^m g(a_i)$ attains the minimum (resp. maximum) at $\bm{a}$.
\begin{theorem}
\label{theorem4}
Assume that there is a number $c \in \{2,\dots,n\}$ such that
\begin{align}
\label{concave}
2g(a) & > g(a-1)+g(a+1) \text{ for all }a \in \{2,\dots,c-1\}\,,\\
\label{lin}
2g(a) & \le g(a-1)+g(a+1) \text{ if }a=c\,,\\
\label{convex}
2g(a) & < g(a-1)+g(a+1) \text{ for all }a \in \{c+1,\dots,n-1\}\,.
\end{align}
\renewcommand{\labelenumi}{\alph{enumi})}~
\begin{enumerate}
\item 
If $\bm{a}^* \in \Pnm$ is a $g$-minimal partition then it has the following properties:

$\bm{a}^*$ does not contain two items from $\{c-1,\dots,n\}$ with absolute difference $\ge 2$,
with the exception that items $c-1$ and $c+1$ may exist if (\ref{lin}) is satisfied with equality.\\
$\bm{a}^*$ does not contain two items from $\{2,\dots,c-1\}$.

\item 
If $\bm{a}^{**} \in \Pnm$ is a $g$-maximal partition then it has the following properties:

$\bm{a}^{**}$ does not contain two items from $\{1,\dots,c\}$ with absolute difference $\ge 2$.\\
$\bm{a}^{**}$ does not contain two items from $\{c,\dots,n-1\}$,
with the exception that two items $c$ may exist if (\ref{lin}) is satisfied with equality.
\end{enumerate}
\end{theorem}
\proof We prove only part a) because part b) can be proved analogously.
Let briefly $g(\bm{a})=\sum_{i=1}^m g(a_i)$ and let $\bm{a}^*$ be a $g$-minimal partition.
Assume that $\bm{a}^*$ contains two items $a_i,a_j$ from $\{c-1,\dots,n\}$ with $a_j-a_i \ge 2$ (where $a_i=c-1$ and $a_j=c+1$ are not considered if (\ref{lin}) is satisfied with equality).
Then let $\bm{a}'$ be the partition that can be obtained from $\bm{a}^*$ by incrementing $a_i$ and decrementing $a_j$. Then
\[
g(\bm{a}')-g(\bm{a}^*)=(g(a_j-1)-g(a_j))-(g(a_i)-g(a_i+1)) < 0
\]
since by (\ref{convex}) for all $a \in \{c+1,\dots,n\}$
\[
g(a-1)-g(a) > g(a)-g(a+1) > g(a+1)-g(a+2) > \cdots
\]
Hence $\bm{a}^*$ cannot be $g$-minimal, a contradiction.

Assume that $\bm{a}^*$ contains two items $a_i,a_j$ from $\{2,\dots,c-1\}$ with $a_j \ge a_i$.
Let $\bm{a}'$ be the partition that can be obtained from $\bm{a}^*$ by decrementing $a_i$ and incrementing $a_j$. Then
\[
g(\bm{a}')-g(\bm{a}^*)=(g(a_j+1)-g(a_j))-(g(a_i)-g(a_i-1)) < 0
\]
since by (\ref{concave}) for all $a \in \{2,\dots,c-1\}$
\[
g(a)-g(a-1) > g(a+1)-g(a) > g(a+2)-g(a+1) > \cdots > g(c)-g(c-1).
\]
Hence $\bm{a}^*$ cannot be $g$-minimal, a contradiction.

\qed

In the  following we write for a $k$-tuple $(a_{i_1},\dots,a_{i_k})$ of items from $\bm{a}$ briefly
$g(a_{i_1},\dots,a_{i_k}) = \sum_{l=1}^k g(a_{i_l})$. In order to obtain contradictions we replace such a $k$-tuple by some $k$-tuple $(a_{i_1}',\dots,a_{i_k}')$
of items such that $a_{i_1}+\dots+a_{i_k}=a_{i_1}'+\dots+a_{i_k}'$ and $g(a_{i_1},\dots,a_{i_k}) <$ (resp. $>$) $g(a_{i_1}',\dots,a_{i_k}')$.

\section{Proof of Theorem \ref{theorem1}}
We have
\[
F_{w_1}(a)=\sum_{k=0}^n \left(\binom{n-a}{k} + \binom{n-a}{k-a}\right)=2^{n-a+1}.
\]
Clearly we may ignore the factor $2^{n+1}$ and thus we work with the function
\[
g(a)=\frac{a}{2^a} \quad \left(=\frac{1}{2^{n+1}} a F_{w_1}(a)\right)\,.
\]
It is easy to check that (\ref{concave}), (\ref{lin}) (with equality) and  (\ref{convex}) of Theorem \ref{theorem4} are satisfied with $c=3$.
Moreover, note that
\begin{equation}
\label{mon}
g(1)=g(2) > g(3) > g(4) > \cdots
\end{equation}

Part a): Let $\bm{a}^*$ be a 1-minimal partition. 
In addition, we choose $\bm{a}^*$ in such a way that it has a minimum number of items 4 with respect to all 1-minimal partitions.
Then, by Theorem \ref{theorem4} a), $\bm{a}^*$ does not contain two items 2. Moreover, by the special choice of $\bm{a}^*$ it does not contain a pair $(2,4)$ because it could be replaced by
$(3,3)$ which contradicts the 1-minimality and the minimality of the number of items 4. Thus, by Theorem \ref{theorem4} a), $\bm{a}^*$ has the form 
$1^i b^j (b+1)^k$ with some integer $b \ge 2$, i.e.,
\begin{equation}
\label{mn}
i+j+k=m \text{ and } i+bj + (b+1) k = n\,.
\end{equation}
We may assume that $k \ge 1$ because a partition $1^ib^j (b+1)^0$ can also be considered as a partition $1^i (b-1)^0 ((b-1)+1)^j$.

Since the reasoning in the proof of Theorem \ref{theorem2} will be analogous we point out that we use in the following the inequalities
\begin{align}
\label{23,14}
g(2,3) & > g(1,4)\,,\\
\label{33,15}
g(3,3) & > g(1,5)\,,\\
\label{34,16}
g(3,4) & > g(1,6)\,,\\
\label{12k-1,kk}
g(1,2k-1) & > g(k,k)\text{ for all } k \ge 4\,,\\
\label{12k,kk+1}
g(1,2k) & > g(k,k+1)\text{ for all } k \ge 4\,,\\
\label{166,445}
g(1,6,6) & > g(4,4,5)\,,\\
\label{444,156}
g(4,4,4) & \ge g(1,5,6)\text{ (in fact $g(4,4,4) = g(1,5,6)$)} \,.
\end{align}
These inequalities can be verified for $g(a)=\frac{a}{2^a}$ by elementary computations.
Concerning (\ref{12k-1,kk}) (and analogously (\ref{12k,kk+1})) note that
$g(1,2k-1)>g(1) =\frac{1}{2} = 2g(4) \ge g(k,k)$).

{\bf Case 1}. $n-m \le 5$. We have to show that $j+k \le 1$. Assume that $j+k \ge 2$. By (\ref{mn}),
\[
5 \ge n-m=(b-1)j+bk \ge (b-1)(j+k) \ge (b-1)2
\]
and hence $b \le 3$ and, moreover, $k=1$ if $b=3$.
Thus $\bm{a}^*$ contains a pair of items of the form $(2,3)$ or $(3,3)$ or $(3,4)$. But
such a pair could be replaced by $(1,4)$, $(1,5)$ or $(1,6)$, respectively,
which contradicts the 1-minimality of $\bm{a}^*$ by (\ref{23,14}), (\ref{33,15}) and (\ref{34,16}), respectively.

{\bf Case 2}. $n-m > 5$ and $n > 5m$. We have to show that $i=0$.

 Assume that $i>0$, i.e., $\bm{a}^*$ contains an item 1. The partition $\bm{a}^*$ cannot contain an item $d \ge 7$ 
because the pair $(1,d)$ could be replaced by $(\frac{d+1}{2},\frac{d+1}{2})$ if $d$ is odd and by $(\frac{d}{2},\frac{d}{2}+1)$
if $d$ is even which contradicts the 1-minimality of $\bm{a}^*$ by (\ref{12k-1,kk}) and (\ref{12k,kk+1}).
Thus $b\le 5$.

Assume that $\bm{a}^*$ contains at most one item 6. Then by (\ref{mn})
\[
5m-n=4i+(5-b)j+(5-(b+1))k \ge 4+0-1>0,
\]
a contradiction to $n > 5m$.
Thus $\bm{a}^*$ contains an item 1 and two items 6. But $(1,6,6)$ can be replaced by $(4,4,5)$ which contradicts the 1-minimality of $\bm{a}^*$ by (\ref{166,445}).

{\bf Case 3}. $n-m > 5$ and $n \le 5m$. 
If $n-m=6$, then, similarly to Case 1, $\bm{a}^*$ contains a pair of items $(2,3)$ or $(3,3)$ or $(3,4)$ or $(4,4)$. Since the first three pairs yield a contradiction $\bm{a}^*$ has the form $1^{m-2} 4^2$ as asserted. 

Thus let $n-m>6$.

First note that $\bm{a}^*$ cannot contain three items 4 because the triple $(4,4,4)$ could be replaced by $(1,5,6)$ which contradicts the 1-minimality and the minimality of the number of items 4 of $\bm{a}^*$ by (\ref{444,156}).
Now we show that 
\begin{equation}
\label{451}
b=4 \text{ or } (b=5 \text{ and } k=1).
\end{equation}

In order to show that $b \ge 4$ assume that $b \le 3$. Then $j+k > 2$ because otherwise $6 < n-m = (b-1)j+bk \le b(j+k) \le 6$. Thus $\bm{a}^*$ contains a pair of items of the form $(2,3)$ or $(3,3)$ or $(3,4)$ which leads as in Case 1 to a contradiction.

Assume that $b \ge 6$. Then $\bm{a}^*$ must contain an item 1 in view of $n \le 5m$ and an item $d=(b+1) \ge 7$. This leads as in Case 2 to a contradiction.

Now let $b=5$. If $i=0$ then in view of $n \le 5m$ necessarily $n=5m$ and hence $\bm{a}^*$ has the form $1^04^05^m$, i.e., $b=4$, a contradiction.
Thus $i>0$, i.e., $\bm{a}^*$ contains an item 1. Then $\bm{a}^*$ cannot contain two items 6 because otherwise $(1,6,6)$ could be replaced by $(4,4,5)$ which contradicts the 1-minimality of $\bm{a}^*$ by (\ref{166,445}). Hence $k= 1$ and (\ref{451}) is proved.


Thus $\bm{a}^*$ has necessarily one of the following forms: $1^i5^k$, $1^i4^15^k$, $1^i4^25^k$, $1^i5^j6^1$.

From (\ref{mn}) it follows that $i = m - \frac{n-m+j}{b}$ and $k=\frac{n-m-(b-1)j}{b}$. Consequently:

The form $1^i5^k$ requires $m-n \equiv 0(\hspace{-1ex}\mod 4)$ and is more precisely given by $1^{m-\frac{n-m}{4}} 5^{\frac{n-m}{4}}$.

The form $1^i4^15^k$ requires $m-n \equiv 1(\hspace{-1ex}\mod 4)$ and is more precisely given by $1^{m-\frac{n-m+1}{4}} 4^1 5^{\frac{n-m-3}{4}}$.

The form $1^i4^25^k$ requires $m-n \equiv 2(\hspace{-1ex}\mod 4)$ and is more precisely given by $1^{m-\frac{n-m+2}{4}} 4^2 5^{\frac{n-m-6}{4}}$.

The form $1^i5^j6^1$ has $1=k=\frac{n-m-(b-1)j}{b}$ which leads to the precise form $1^{m-\frac{n-m-1}{4}} 5^{\frac{n-m-5}{4}} 6^1$ and requires 
$m-n \equiv 3(\hspace{-1ex}\mod 4)$.

Thus part a) is proved.

Part b): Let $\bm{a}^{**}$ be a 1-maximal partition.
In addition, we choose $\bm{a}^{**}$ in such a way that it has a minimum number of items 3 with respect to all 1-maximal partitions.
Then, by Theorem \ref{theorem4} b), $\bm{a}^{**}$ has the form 
$1^i 2^j b^k$ with some integer $b \ge 3$ and $k \le 1$.

{\bf Case 1}. $n>2m$. Then $k=1$. We have to show that $i=0$. Assume that $i \ge 1$. By (\ref{mon}), $g(1,b)<g(2,b-1)$ which contradicts the 1-maximality of $\bm{a}^{**}$.

{\bf Case 2}. $n\le 2m$. We have to show that $k=0$. Assume that $k \ge 1$. 
If $i=0$ then $n=2j+bk > 2(j+k)=2m$, a contradiction. Thus $i \ge 1$. As in Case 1, $(1,b)$ can be replaced by $(2,b-1)$, a contradiction.

\qed

\section{Proof of Theorem \ref{theorem2}}
We have
\[
F_{w_2}(a)=\sum_{k=1}^n  \left( \binom{n-a}{k} + \binom{n-a}{k-a}\right)\frac{1}{k}\,.
\]
The main idea is to use an integral representation of $F_{w_2}(a)$ and to prove inequalities ``under the integral''.

Let
\begin{align*}
p_a^+(t)&=\frac{1}{t}\left((1+t)^{n-a}-1\right),\\
p_a^-(t)&=t^{a-1}(1+t)^{n-a}\,.
\end{align*}
The binomial theorem and a simple evaluation of the integrals yield the following representation.

\begin{lemma}
\label{lemma2.1}
We have
\[
F_{w_2}(a)=\int_0^1 (p_a^+(t)+p_a^-(t))\,dt\,.
\]
\end{lemma}

Consequently, for the proof of Theorem \ref{theorem2}, we work with the function
\[
g(a)=\int_0^1 a (p_a^+(t)+p_a^-(t))\,dt\,.
\] 

In order to prove inequalities, we need the following estimation.

\begin{lemma}
\label{lemma2.3} 
Let $\varphi(t)$ and $\varrho(t)$ be continuous real functions on $[0,1]$. Let $\varphi(t)$ be positive and increasing and suppose that there is some root $x_0 \in (0,1)$ of $\varrho(t)$ such that $\varrho(t)< 0$ for all $t \in [0,x_0)$ and $\varrho(t)> 0$ for all $t \in (x_0,1]$.
If
\begin{equation}
\label{rho}
\int_0^1 \varrho(t)\,dt > 0\quad (\text{resp.} \int_0^1 \varrho(t)\,dt \ge 0)
\end{equation}
then
\begin{equation}
\label{phirho}
\int_0^1 \varphi(t)\varrho(t)\,dt > 0 \quad (\text{resp.} \int_0^1 \varphi(t)\varrho(t)\,dt \ge 0)\,.
\end{equation}
\end{lemma}
\proof
We have
\begin{align*}
\int_0^1 \varphi(t)\varrho(t)\,dt&=\int_0^{x_0} \varphi(t)\varrho(t)\,dt+\int_{x_0}^1 \varphi(t)\varrho(t)\,dt\\
&\ge \int_0^{x_0} \varphi(x_0)\varrho(t)\,dt+\int_{x_0}^1 \varphi(x_0)\varrho(t)\,dt\\
&=\varphi(x_0)\int_0^{1} \varrho(t)\,dt\\
&> 0 \quad (\text{resp.} \ge 0)\,.
\end{align*}
\qed

First we study the monotonicity of $g$.

\begin{lemma}
\label{lemmamon}
Let $n \ge 6$.
We have
\begin{align*}
g(a+1) &> g(a) \text{ if } a=1\,,\\
g(a+1) &< g(a) \text{ if } a \in \{2,\dots,n-1\}\,.
\end{align*}
\end{lemma}
\proof In fact, the functions $p_a^+, p_a^-,g$ depend also on $n$, hence we write here more detailed $p_{a,n}^+, p_{a,n}^-, g(a,n)$.
Moreover, let
\begin{align*}
p_{a,n}(t)&=a (p_{a,n}^+(t)+p_{a,n}^-(t))\,,\\
q_{a,n}(t)&=p_{a,n+1}(t)-p_{a,n}(t)\,,\\
d_{a,n}(t)&=p_{a+1,n}(t)-p_{a,n}(t)\,.\\
\end{align*}
Note that 
\begin{align*}
g(a,n)&=\int_0^1 p_{a,n}(t)\,dt\,,\\
q_{a,n}&=a(1+t)^{n-a}(1+t^a)\,
\end{align*}
and
\begin{equation}
\label{gdiff}
g(a+1,n)-g(a,n)=\int_0^1 d_{a,n}(t)\,dt\,.
\end{equation}
By definition, 
\begin{align*}
p_{a,a+1}&=a+a t^{a-1}(1+t)\,,\\
p_{a+1,a+1}&=(a+1) t^a\,,
\end{align*}
and hence, for all $t \in [0,1]$,
\[
d_{a,a+1}=-a+t^{a-1}(t-a) < 0\,,
\]
which implies together with (\ref{gdiff}) that
\begin{equation}
\label{init}
g(a+1,a+1)-g(a,a+1) < 0.
\end{equation}
For $n \ge a+1$ we have
\begin{align*}
d_{a,n+1}(t)-d_{a,n}(t)&= \left(p_{a+1,n+1}(t)-p_{a,n+1}(t)\right)-\left(p_{a+1,n}(t)-p_{a,n}(t)\right)\\
&=q_{a+1,n}(t)-q_{a,n}(t)\\
&=(a+1)(1+t)^{n-a-1}(1+t^{a+1})-a(1+t)^{n-a}(1+t^a)\\
&=(1+t)^{n-a-1}\left((1-at)+t^a(t-a)\right)\,.
\end{align*}
For $a=1$ we have for all $t \in (0,1]$
\[
d_{1,n+1}(t)-d_{1,n}(t)=(1+t)^{n-2}(1-t)^2 > 0
\]
and by (\ref{gdiff})
\[
\left(g(2,n+1)-g(1,n+1)\right)-\left(g(2,n)-g(1,n)\right) > 0\,.
\]
Since $g(2,6)-g(1,6)= \frac{29}{60} > 0$ it follows that
$g(2,n) > g(1,n)$ if $n \ge 6$.

For $a \ge 2$ we have 
\[
d_{a,n+1}(t)-d_{a,n}(t) \ge (1+t)^{n-a-1}\left((1-2t)+t^a(t-a)\right)
\]
By Lemma \ref{lemma2.3} applied to $\varphi(t)=(1+t)^{n-a-1}$ and $\varrho(t)=2t-1$ and in view of $t^a(t-a) < 0$ for all $t \in (0,1]$ we obtain easily that
\[
\int_0^1 (1+t)^{n-a-1}\left((1-2t)+t^a(t-a)\right) \,dt < 0\,,
\]
and hence
\[
\left(g(a+1,n+1)-g(a,n+1)\right)-\left(g(a+1,n)-g(a,n)\right) < 0\,.
\]
which together with (\ref{init}) implies that $g(a+1,n) < g(a,n)$ if $a \ge 2$. 

\qed

Let $\bm{a}=(a_1,\dots,a_k)$ and $\bm{\alpha}=(\alpha_1,\dots,\alpha_k)$ be  $k$-tuples of integers.
To prove further inequalities, we have to evaluate sums of the form $\sum_{i=1}^k \alpha_ig(a_i)$. For this purpose, we introduce the polynomials

\begin{align*}
q_{\bm{\alpha},\bm{a}}^+(t)&=\sum_{i=1}^k \alpha_i a_i (1+t)^{a_k-a_i}\,,\\
q_{\bm{\alpha},\bm{a}}^-(t)&=\sum_{i=1}^k \alpha_i a_i t^{a_i-a_1}(1+t)^{a_k-a_i}\,,
\end{align*}
but, for a simpler representation, we omit the indices $\bm{\alpha},\bm{a}$ and write simply $q^+(t)$ and $q^-(t)$.

Note  that
\begin{equation}
q^+(t)=t^{a_k-a_1}q^-(1/t)\,.
\end{equation}

\begin{lemma}
\label{lemma2.2}
Let $\sum_{i=1}^k \alpha_ia_i=0$. Then
\[
\sum_{i=1}^k \alpha_ig(a_i)=\int_0^1 (1+t)^{n-a_k}\left(t^{a_k-a_1-1}q^-(1/t)+ t^{a_1-1}q^-(t)\right)\,dt\,.
\]
\end{lemma}
\proof We have
\begin{align*}
\sum_{i=1}^k \alpha_ia_i(p_{a_i}^+(t)+p_{a_i}^-(t))&=\sum_{i=1}^k \alpha_ia_i \frac{1}{t}\left((1+t)^{n-a_i}-1\right)+\alpha_ia_it^{a_i-1}(1+t)^{n-a_i}\\
&=(1+t)^{n-a_k} \sum_{i=1}^k \alpha_ia_i \frac{1}{t} (1+t)^{a_k-a_i}+\alpha_ia_it^{a_i-1}(1+t)^{a_k-a_i}\\
&=(1+t)^{n-a_k} \left(\frac{1}{t} q^+(t)+t^{a_1-1}q^-(t)\right)\\
&=(1+t)^{n-a_k} \left(t^{a_k-a_1-1}q^-(1/t)+t^{a_1-1}q^-(t)\right)\,.
\end{align*}
\qed

In the following let
\[
r(t)=t^{a_k-a_1-1}q^-(1/t)+t^{a_1-1}q^-(t)
\]
which implies that
\begin{equation}
\label{intest}
\sum_{i=1}^k \alpha_ig(a_i)=\int_0^1 (1+t)^{n-a_k}r(t)\,dt\,.
\end{equation}


\begin{lemma}
\label{lemma>0}~
\renewcommand{\labelenumi}{\alph{enumi})}~
\begin{enumerate}
\item If $r(t)\ge 0$ for all $t \in [0,1]$ and $r(t)$ is not the zero polynomial then $\sum_{i=1}^k \alpha_ig(a_i)>0$ for all $n \ge a_k$.
\item If there is some root $x_0 \in (0,1)$ of $r(t)$ such that $r(t)< 0$ for all $t \in [0,x_0)$ and $r(t)> 0$ for all $t \in (x_0,1]$ and if there is some integer $l$ such that
$\int_0^1 (1+t)^{l}r(t)\,dt>0$ then $\sum_{i=1}^k \alpha_ig(a_i)>0$ for all $n \ge a_k+l$.
\end{enumerate}
\end{lemma}
\proof Part a) is a trivial consequence of (\ref{intest}) and part b) follows from Lemma \ref{lemma2.3} applied to $\varphi(t)=(1+t)^{n-a_k-l}$ and $\varrho(t)=(1+t)^{l}r(t)$.

\qed

\begin{lemma}
\label{lemma2.4} 
Let $n \ge 6$.
We have
\begin{align*}
2g(a) &> g(a-1)+g(a+1) \text{ if } a \in \{2,3\}\,,\\
2g(a) &< g(a-1)+g(a+1) \text{ if } a \in \{4,\dots,n-1\}\,.
\end{align*}
\end{lemma}
\proof We have to estimate $g(a-1)-2g(a)+g(a+1)$ and hence we work with $\bm{a}=(a-1,a,a+1)$ and $\bm{\alpha}=(1,-2,1)$ and apply Lemma \ref{lemma>0}. We have
\begin{align*}
q^-(t)&=(a-1)(1+t)^2-2at(1+t)+(a+1)t^2\,,\\
r(t)&
=t^{a-2}(a-2t-1)+(a-1)t-2\,.
\end{align*}
If $a=2$ then $r(t)=-(t+1)< 0$ for all $t \in [0,1]$, if $a=3$ then $r(t)=-2(t-1)^2<0$ for all $t \in [0,1)$ and hence by Lemma \ref{lemma>0} a) (applied to $-r(t)$)
\[
g(a-1)-2g(a)+g(a+1)< 0\quad \text{ if } a \in \{2,3\}\,.
\]
If $a=4$ then $r(t)=-2t^3+3t^2+3t-2$, $r(t)<0$ for all $t \in [0,\frac{1}{2})$, $r(t)>0$ for all $t \in (\frac{1}{2},1]$, and $\int_0^1 (1+t)^1r(t)\,dt=\frac{7}{20}>0$.
By Lemma \ref{lemma>0} b) (note that $n-(a+1) \ge 6-5=1$)
\[
g(a-1)-2g(a)+g(a+1) > 0\quad \text{ if } a =4\,.
\]
Now let $a \ge 5$. Clearly $a-2t-1 > 0$ for all $t \in [0,1]$ and hence
\[
\int_0^1 r(t)\,dt > \int_0^1 (a-1)t-2 \ge 0\,.
\]
By Lemma \ref{lemma>0} a) 
\[
g(a-1)-2g(a)+g(a+1) > 0\quad \text{ if } a \in \{5,\dots,n-1\}\,.
\] 
\qed

The next lemma summarizes several further inequalities that are needed in the following.
\begin{lemma}
\label{lemmainequ}
Let $n \ge 13$. Then
\begin{align}
\label{2,17,44}
g(1,7) & > g(4,4)\,,\\
\label{2,18,45}
g(1,8) & > g(4,5)\,,\\
\label{2,23,14}
g(2,3) & > g(1,4)\,,\\
\label{2,25,34}
g(2,5) & > g(3,4)\,,\\
\label{2,33,15}
g(3,3) & > g(1,5)\,,\\
\label{2,34,16}
g(3,4) & > g(1,6)\,,\\
\label{2,166,445}
g(1,6,6) & > g(4,4,5)\,,\\
\label{2,225,333}
g(2,2,5) & > g(3,3,3)\,,\\
\label{2,444,156}
g(4,4,4) & > g(1,5,6) \,.
\end{align}
\end{lemma}
\proof For the proof, we again apply Lemma \ref{lemma>0}, analogously to the proof of Lemma \ref{lemma2.4}.
For the sake of brevity we present the corresponding vectors $\bm{a}$ and $\bm{\alpha}$ in form of a table.

\begin{tabular}{c|c|c}
inequality&$\bm{a}$&$\bm{\alpha}$\\\hline
$g(1,7)  > g(4,4)$ & $(1,4,7)$ & $(1,-2,1)$\\
$g(1,8)  > g(4,5)$ & $(1,4,5,8)$ & $(1,-1,-1,1)$\\
$g(2,3)  > g(1,4)$ & $(1,2,3,4)$ & $(-1,1,1,-1)$\\
$g(2,5)  > g(3,4)$ & $(2,3,4,5)$ & $(1,-1,-1,1)$\\
$g(3,3)  > g(1,5)$ & $(1,3,5)$ & $(-1,2,1)$\\
$g(3,4)  > g(1,6)$ & $(1,3,4,6)$ & $(-1,1,1,-1)$\\
$g(1,6,6)  > g(4,4,5)$ & $(1,4,5,6)$ & $(1,-2,-1,2)$\\
$g(2,2,5)  > g(3,3,3)$ & $(2,3,5)$ & $(2,-3,1)$\\
$g(4,4,4)  > g(1,5,6)$ & $(1,4,5,6)$ & $(-1,3,-1,-1)$
\end{tabular}

The corresponding terms $r(t)$ (if Lemma \ref{lemma>0} a) can be applied) and additionally the values $l, \int_0^1 (1+t)^lr(t)\,dt$ and the lower bound $a_k+l$ for $n$ (if Lemma \ref{lemma>0} b) can be applied) can be computed and checked using a standard computer algebra system. They are given in the next table in the corresponding order.

\begin{tabular}{c|c|c|c}
$r(t)$&$l$&$\int_0^1 (1+t)^lr(t)\,dt$&$a_k+l$\\\hline
$-17t^5-3t^4+27t^3+27t^2-3t-17$ &4& $1639/360$&$11$\\
$-23t^6-11t^5+35t^4+62t^3+35t^2-11t-23$&$0$&$23/42$&$8$\\
$3t^2-2t+1/3\ge 0$&&&\\
$-4t^3+5t^2+5t-4$&$3$&$37/70$&$8$\\
$7t^3-4t^2-4t+7\ge 0$&&&\\
$11t^4-2t^3-14t^2-2t+11\ge 0$&&&\\
$-15t^4+7t^3+20t^2+7t-15$&$7$&$247/660$&$13$\\
$-6t^3+7t^2+7t-6$&$5$&$443/168$&$10$\\
$13t^4-3t^3-20t^2-3t+13\ge 0$&&&
\end{tabular}

\qed

\begin{lemma}
\label{lemmakk}
Let $n \ge 11$. We have for all $k \ge 4$ with $2k \le n$ resp. $2k <n$
\begin{align}
\label{2,12k-1,kk}
g(1,2k-1) & > g(k,k)\,,\\
\label{2,12k,kk+1}
g(1,2k) & > g(k,k+1)\,.
\end{align}
\end{lemma}
\proof Both inequalities can be proved by induction on $k$ using (\ref{2,17,44}) and (\ref{2,18,45}) as the induction basis. As an example we show the step from $k$ to $k+1$ for (\ref{2,12k-1,kk}).
We have for $k \ge 4$ by Lemma \ref{lemma2.4}
\begin{align*}
g(2k+1)-g(2k-1)&=\left(g(2k+1)-g(2k)\right) + \left(g(2k)-g(2k-1)\right)\\
&\ge \left(g(k+1)-g(k)\right)+\left(g(k+1)-g(k)\right)
\end{align*}
and hence, by the induction hypothesis,
\[
g(1,2k+1)  - g(k+1,k+1) \ge g(1,2k-1)  - g(k,k) >0\,.
\]
\qed

\begin{lemma}
\label{lemma2.5}
Let $n \ge 8$. Then for all $a$ with $4\le a \le n-3$
\begin{equation}
\label{2,2a+1,3a}
g(2)+g(a+1)> g(3)+g(a)\,.
\end{equation}
\end{lemma}
\proof
This inequality can be proved by induction on $a$ using (\ref{2,25,34}) as the induction basis.
For the induction step from $a$ to $a+1$ note that by Lemma \ref{lemma2.4}, $g(a+2)-g(a+1) > g(a+1)-g(a)$. Using the induction hypothesis we obtain
\[
g(a+2)-g(a+1)+g(2)-g(3) > g(a+1)-g(a)+g(2)-g(3) > 0.
\]
\qed

{\it Proof of Theorem \ref{theorem2}}.

Part a): The case $n=m+1$ is trivial, thus let $n > m+1$. Recall that $n \ge 13$.
Let $\bm{a}^*$ be a 2-minimal partition. Assume that it contains an item 2. In view of $n > m+1$ it contains also an item $a \ge 2$. But the pair $(2,a)$ can be replaced by $(1,a+1)$, a contradiction to the 2-minimality of $\bm{a}^*$ by Lemma \ref{lemmamon}. 
The conditions (\ref{concave}), (\ref{lin}) (with strict inequality) and (\ref{convex}) of Theorem \ref{theorem4} are satisfied for $c=4$ by Lemma \ref{lemma2.4}.
Thus by Theorem \ref{theorem4}, $\bm{a}^*$ has the form $1^ib^j (b+1)^k$ with some $b \ge 3$.

Now we are in the same situation as in the proof of part a) of Theorem \ref{theorem1} (with the stronger inequality $b \ge 3$ instead of $b \ge 2$).
The inequalities (\ref{23,14})--(\ref{444,156}) are satisfied for the actual function $g(a)=a F_{w_2}(a)$ by Lemmas \ref{lemmainequ} and \ref{lemmakk}.
Thus the Cases 1--3 can be discussed exactly in the same way as in the proof of part a) of Theorem \ref{theorem1}.
Note that we do not need the additional assumption that $\bm{a}^*$ contains a minimal number of items 4 because the inequality (\ref{444,156}) is here a strict inequality, see (\ref{2,444,156}).

Part b):
Let $\bm{a}^{**}$ be a 2-maximal partition. It cannot contain a pair $(3,a)$ with $a \ge 4$ since it could be replaced by $(2,a+1)$ which contradicts the 2-maximality of $\bm{a}^{**}$ by Lemma \ref{lemma2.5}.
Moreover, $\bm{a}^{**}$ cannot contain a pair $(1,a)$ with $a \ge 3$ because it could be replaced by $(2,a-1)$ (Lemma \ref{lemmamon}).

As for part a), we may apply Theorem \ref{theorem4} (but here part b)) with $c=4$ and strict inequality in (\ref{lin}). Thus $\bm{a}^{**}$ does not contain two items from $\{1,2,3,4\}$ with absolute difference $\ge 2$ and contains at most one item from $\{4,\ldots,n-1\}$.


{\bf Case 1.1} $n>2m+2$. Assume that there is no item $a \ge 4$. Then there must be at least three items 3. But the triple $(3,3,3)$ could be replaced by $(2,2,5)$
using (\ref{2,225,333}).
Thus there is an item $a \ge 4$ which immediately leads to $2^{m-1}(n-2m+2)^1$.

{\bf Case 1.2} $n=2m+2$. If there is an item $a \ge 5$ then there is no item 1 and hence $n \ge 2(m-1)+a >2m+2$, a contradiction.
Thus the only possibility is $2^{m-2} 3^2$.

{\bf Case 1.3} $n=2m+1$. If there is an an item $a \ge 4$ then there is no item 1 and hence $n \ge 2(m-1)+a >2m+1$, a contradiction.
Thus the only possibility is $2^{m-1} 3^1$.

{\bf Case 2} $n \le 2m$. If $\bm{a}^*$ contains an item $a \ge 3$ then it must also contain an item 1, a contradiction.
Thus $\bm{a}^*$ contains only items 1 and 2 which yields the assertion.

\qed

\section{Proof of Theorem \ref{theorem3}}
We have
\[
F_{w_3}(a)=\sum_{k=1}^n \left( \binom{n-a}{k} + \binom{n-a}{k-a}\right)\frac{1}{\binom{n}{k}k}\,.
\]
We need a much simpler representation of $F_{w_3}(a)$ and apply the following lemma from \cite{Thu07}.

\begin{lemma}
\label{lemmathu}
If $a,b,c$ are integers such that $a \ge 0, b > 0, c \ge a+b$ then
\[
\sum_{j=0}^a \binom{a}{j} \frac{1}{(b+j)\binom{c}{b+j}}= \frac{1}{b\binom{c-a}{b}}\,.
\]
\end{lemma}
In fact, we need this lemma only in the case $c=a+b$ and in the case $b=1$. In these cases the lemma can be proved in a straightforward way.

\begin{lemma}
\label{lemma2}We have
\[
\sum_{k=1}^n \binom{n-a}{k-a}\frac{1}{\binom{n}{k}k} =\frac{1}{a}\,.
\]
\end{lemma}
\proof

We have
\[
\sum_{k=1}^n \binom{n-a}{k-a}\frac{1}{\binom{n}{k}k}=\sum_{j=0}^{n-a} \binom{n-a}{j}\frac{1}{(a+j)\binom{n}{a+j}}
\]
and may apply Lemma \ref{lemmathu} with $a:=n-a, b:=a, c:=n$.

\qed

%
%
%
%

\begin{lemma}
\label{lemma5}We have
\[
\sum_{k=1}^n \binom{n-a}{k}\frac{1}{\binom{n}{k}k} =\sum_{k=a+1}^n\frac{1}{k}\,.
\]
\end{lemma}
\proof

We have 
\begin{align*}
\sum_{k=1}^n \binom{n-a}{k}\frac{1}{\binom{n}{k}k}&=\sum_{k=1}^n \left(\sum_{l=0}^{n-a-1}\binom{l}{k-1}\right)\frac{1}{\binom{n}{k}k}\\
&=\sum_{l=0}^{n-a-1} \left(\sum_{k=1}^n \binom{l}{k-1}\frac{1}{\binom{n}{k}k}\right)\\
&=\sum_{l=0}^{n-a-1} \left(\sum_{j=0}^l \binom{l}{j}\frac{1}{(1+j)\binom{n}{1+j}}\right)\\
&=\sum_{l=0}^{n-a-1} \frac{1}{n-l}\\
&=\sum_{k=a+1}^n \frac{1}{k}\,.
\end{align*}
Here the second to last identity follows from Lemma \ref{lemmathu} with $a:=l, b:=1, c:=n$.

\qed

%

Now we are ready to present the simpler form of $F_{w_3}(a)$ which follows immediately from  Lemmas \ref{lemma2} and \ref{lemma5}.
\begin{theorem}
\label{theorem5}
We have
\[
F_{w_3}(a)=\frac{1}{a}+\sum_{k=a+1}^{n} \frac{1}{k}\,.
\]
\end{theorem}

Consequently, for the proof of Theorem \ref{theorem3}, we work with the function
\[
g(a)=1+a \sum_{k=a+1}^{n} \frac{1}{k} \quad \left(=a F_{w_3}(a)\right)\,.
\] 
We have for all $a \in \{2,\dots,n-1\}$
\[
g(a-1)+g(a+1)-2g(a)=(a-1)\left(\frac{1}{a}+\frac{1}{a+1}\right)-\frac{2a}{a+1} = -\frac{1}{a} < 0\,,
\]
and hence we may apply Theorem \ref{theorem4} with $c=n$ and this theorem immediately yields the statements in Theorem  \ref{theorem3}.

\qed

\section{Concluding remarks}

We proved Theorem \ref{theorem2} only for $n \ge 13$. For smaller $n$, some of the used inequalities are not true. which leads to some differences.
But for $3 \le n \le 12$, the values $g(a)$, $a \in \{1,\dots,n-1\}$, can be concretely determined and optimal partitions can be determined by complete search or using a slight modification of the algorithm from \cite{ERS14}.
We omit the details.

With analogous methods, the average cardinality of the $\A$-upper approximation and of the $\A$-lower approximation, respectively, can be optimized. In the case $j=2$ (and in a trivial way also for $j=3$) there are some small differences concerning the used inequalities which leads to some small differences concerning the optimal partitions. We again omit the details.

An inspection of the proof of Lemma \ref{lemma1} and the application of Lemma \ref{lemma2} lead to the identity
\[
\sum_{X \subseteq [n]} w_3(|X|) |\A^-(X)| = \sum_{i=1}^m a_i \sum_{k=0}^n \binom{n-a_i}{k-a_i} w_3(k)=\sum_{i=1}^m a_i \frac{1}{a_i} = m\,.
\]
Consequently,
\[
\sum_{\emptyset \neq X \subseteq [n]} \frac{|\A^-(X)|}{\binom{n}{|X|}|X|}=m\,.
\]
This is an interesting analogon to the well-known AZ-identity \cite{AZ90} (for a textbook see \cite{Eng97}):
\begin{theorem}[AZ-identity]
Let $\F$ be a family of nonempty subsets of $[n]$. For $X \subseteq [n]$ let 
\[
\F(X)=\bigcap_{A \in \F: A \subseteq X} A
\]
(here the intersection is defined to be empty if there is no $A \in \F$ with $A \subseteq X$). Then
\[
\sum_{\emptyset \neq X \subseteq [n]} \frac{|\F(X)|}{\binom{n}{|X|}|X|}=1\,.
\]
\end{theorem}

\bibliographystyle{elsarticle-num}

\end{document}